\documentclass[12pt]{amsart}
\usepackage{amssymb,amsmath}
\usepackage[dvips]{graphicx}
\usepackage{pslatex}

\theoremstyle{plain} \newtheorem{thm}{Theorem}
 \newtheorem{prop}[thm]{Proposition}

\newtheorem*{namedtheorem}{\theoremname}
\newcommand{\theoremname}{testing}
\newenvironment{named}[1]{\renewcommand{\theoremname}{#1}\begin{namedtheorem}}{\end{namedtheorem}}

\theoremstyle{definition} \newtheorem{defn}[thm]{Definition}

\theoremstyle{remark}

\begin{document}

\title{Factorial growth rates for the number of hyperbolic 3-manifolds of a given volume}

\author{Christian Millichap} 
\address{Department of Mathematics, Temple University\\ Philadelphia, PA 19122}
\email{Christian.Millichap@gmail.com}

\begin{abstract}
The work of J\o rgensen and Thurston shows that there is a finite number $N(v)$ of orientable hyperbolic $3$-manifolds with any given volume $v$.  In this paper, we construct examples showing that  the number of hyperbolic knot complements with a given volume $v$ can grow at least factorially fast with $v$.  A similar statement holds for closed hyperbolic $3$-manifolds, obtained via Dehn surgery.  Furthermore, we give explicit estimates for lower bounds of $N(v)$ in terms of $v$ for these examples.  These results improve upon the work of Hodgson and Masai, which describes examples that grow exponentially fast with $v$.  Our constructions rely on performing volume preserving mutations along Conway spheres and on the classification of Montesinos knots.  
\end{abstract}

\maketitle
\section{Introduction}
\label{sec:intro}

Over the past thirty-five years, the set of volumes of orientable, complete, hyperbolic $3$-manifolds (henceforth, referred to as hyperbolic $3$-manifolds) has been studied extensively.  J\o rgensen and Thurston \cite{Th} have shown that the set of volumes is in fact a closed, non-discrete, well-ordered subset of $\mathbb{R}_{>0}$ of order type $\omega^{\omega}$, and that there is a finite number $N(v)$ of hyperbolic 3-manifolds with any given volume $v$.  Even though this volume function is finite to one, Wielenberg \cite{Wi} has proved that there does not exist a uniform upper bound for the number of hyperbolic $3$-manifolds of a given volume.  In this paper, we construct examples showing that $N(v)$ grows at least factorially fast with $v$ and can take values larger than $v^{\left(\frac{v}{8}\right)}$ for $ v \gg 0$. We also provide evidence supporting the claim that a factorial growth rate is actually the best possible.  

Here, we examine how the volume function behaves on hyperbolic knot complements and closed hyperbolic $3$-manifolds.  The Callahan--Hildebrand--Weeks census of cusped hyperbolic $3$-manifolds \cite{CHW} and the Hodgson--Weeks census of closed hyperbolic $3$-manifolds \cite{HW} provide frequencies of volumes of hyperbolic $3$-manifolds from their respective classes.  While these results are only lower bounds on $N(v)$, they do suggest that often there are only a few manifolds with any given volume.  In fact, Hodgson and Masai prove that there exists an infinite sequence of closed hyperbolic $3$-manifolds determined by their volumes, and similarly, an infinite sequence of $1$-cusped hyperbolic $3$-manifolds determined by their volumes \cite{HM}.      

Another notable feature in the frequencies given by the censuses is that there are particular values of $v$ for which $N(v)$ is very large.  This raises the question, just how fast does $N(v)$ grow with respect to $v$?  Recently, Hodgson and Masai \cite{HM} and Chesebro and DeBlois  \cite{CD} have each constructed examples showing that the number of hyperbolic link complements with a given volume, $N_{L}(v)$, grows at least exponentially fast with $v$.  Hodgson and Masai's constructions rely on volume preserving mutations along thrice punctured spheres, while Chesebro and DeBlois use volume preserving mutations along Conway spheres.  

At the end of \cite{HM}, Hodgson and Masai pose a number of open questions related to the volume function and $N(v)$, one of which asks if there are sequences $v_{n} \rightarrow \infty$ such that the growth rate of $N_{L}(v_{n})$ with respect to $v_{n}$ is faster than exponential.  Faster than exponential growth rates for $N(v)$ have been established for a few other special classes of hyperbolic $3$-manifolds.  For instance, Frigerio--Martelli--Petronio \cite{FMP} have constucted hyperbolic $3$-manifolds with totally geodesic boundary giving sequences $v_{n} \rightarrow \infty$ such that the number $N_{G}(v_{n})$ grows at least as fast as $v_{n}^{cv_{n}}$ for some $c > 0$, where $N_{G}(v)$ denotes the number of hyperbolic $3$-manifolds with totally geodesic boundary and volume $v$.  Belolipetsky--Gelander--Lubotzky--Shalev \cite{BGLS} have shown that the number of arithmetic hyperbolic $3$-manifolds of volume at most $v$ is roughly $v^{c'v}$ for $c' > 0$ and $v \gg 0$.  This result implies that for $v \gg 0$, $\max_{v_{i}\leq v} N_{A}(v_{i}) \approx v^{c'v}$, where $N_{A}(v)$ denotes the number of arithmetic hyperbolic $3$-manifolds of volume $v$.  Though both of these papers show faster than exponential growth rates for $N(v)$, they do not give sharp lower or upper bounds of $N(v)$ in terms of $v$. Also, Lubotzky and Thurston explained how to find factorial growth rates for $N(v)$ by estimating the number of covering spaces of a fixed non-arithmetic hyperbolic $3$-manifold whose fundamental group surjects onto a free group of rank $2$. In principle, explicit estimates could be obtained for lower bounds of $N(v)$ in terms of $v$ using this technique; see \cite{BGLM} and \cite{Ca} for details.  

In this paper, we provide an explicit answer to Hodgson and Masai's question by showing that both the number of knot complements with a given volume and the number of closed $3$-manifolds with a given volume each grow at least factorially fast with volume.  In fact, we show that there are infinitely many sequences of volumes for closed $3$-manifolds that exhibit this growth rate.  On top of this, we are able to give explicit lower bounds for growth rates of $N(v)$ in terms of $v$.  This substantial improvement is highlighted in our main theorem:

\begin{thm}
\label{thm:mainthm}
There exist sequences of volumes, $v_{n} , x_{n} \rightarrow \infty$ such that 
\begin{center}
$N_{K}(v_{n}) \geq \left(v_{n} \right)^{\left(\frac{v_{n}}{8} \right)}$ and $N_{Cl}(x_{n}) \geq \left(x_{n} \right)^{\left(\frac{x_{n}}{8} \right)}$
\end{center}
for all $n \gg 0$. Here, $N_{K}(v)$ is the number of knot complements with volume $v$ and $N_{Cl}(x)$ is the number of closed hyperbolic $3$-manifolds with volume $x$.  
\end{thm}

Theorem \ref{thm:mainthm} is a direct consequence of the following two theorems, which use our constructions of hyperbolic knot complements that share the same volume to show how fast $N_{K}(v)$ and $N_{Cl}(v)$ can grow.

\begin{thm}
\label{thm:knotvol}
There exists a sequence of volumes $\left\{v_{n}\right\}_{n=2}^{\infty}$ such that $N_{K}(v_{n}) \geq \frac{(2n)!}{2}$.
 
Furthermore, $\left(\frac{2n-1}{2}\right)v_{\mathrm{oct}} \leq v_{n}  \leq \left(4n+2\right)v_{\mathrm{oct}}$, where $v_{\mathrm{oct}} \left(\approx 3.6638\right)$ is the volume of a regular ideal octahedron.
\end{thm}

\begin{thm}
\label{thm:closedvol}
For each $n \geq 2$ there exists a convergent sequence $x_{n_{i}} \rightarrow v_{n}$ such that $N_{Cl}(x_{n_{i}}) \geq \frac{(2n)!}{2}$.  Furthermore, $x_{n_{i}} \leq \left(4n+2\right)v_{\mathrm{oct}}$, for all $i$.  
\end{thm}

In Theorem \ref{thm:knotvol}, we construct sequences of volumes $v_{n} \rightarrow \infty$ such that the number of knot complements $N_{K}(v_{n})$ grows at least factorially fast with $v_{n}$.  Similar to Chesebro and DeBlois, we use mutations along Conway spheres as a simple way to construct hyperbolic $3$-manifolds with the same volume.  We are able to distinguish these mutated knots from one another via the classification of Montesinos knots.  Then we perform particular Dehn surgeries on these mutated knots to obtain non-homeomorphic closed hyperbolic $3$-manifolds with the same volume, giving us Theorem \ref{thm:closedvol}. 

Now, let us outline the rest of this paper.  In Section \ref{sec:TandK}, we recall some fundamental results about rational tangles and how they can be used to construct Montesinos knots.  In particular, rational tangles can be fully described by their representative fractions and Montesinos links can be described by a finite sequence of rational tangles, which makes it very simple to describe such links.  Furthermore, using Montesinos links in our constructions makes it easy to distinguish two mutated links from one another via the classification of Montesinos links; see \cite{Bo}, \cite{BS}, and \cite{BZ}.  In section \ref{sec:Csphere}, we review Ruberman's work \cite{Ru} on performing mutations along Conway spheres.  This will serve as our tool for creating new hyperbolic knots that have the same volume.  In section \ref{sec:DSurgery}, we go over some basic facts about performing Dehn surgery on cusped hyperbolic $3$-manifolds, and once again highlight a theorem of Ruberman showing that for particular mutations on knots, corresponding Dehn surgeries preserve volume.  This will be our main tool in constructing closed hyperbolc $3$-manifolds that have the same volume.  Section \ref{sec:Constructions} gives the explicit details for constructing our Montesinos links, which we plan to mutate.  Sections \ref{sec:KnotVol} and \ref{sec:ClosedVol} give the proofs of Theorems \ref{thm:knotvol} and \ref{thm:closedvol}, respectively.  Finally, Section \ref{sec:Final} gives the proof of Theorem \ref{thm:mainthm} and highlights some open problems related to the growth rate of $N(v)$ with respect to $v$.  

We are grateful to David Futer for his help and guidance with this project.

\section{Tangles and knots}
\label{sec:TandK}

Throughout this paper we will be working with \emph{rational tangles} and knots constructed by performing certain basic operations on rational tangles.  The theory of tangles was first studied by John Conway \cite{Co} in his work on knot tabulation.  For more background on rational tangles, see \cite{KL}.      

\begin{defn}[$2$-tangle] 
Let $t$ be a pair of unoriented arcs properly embedded in a $3$-ball $B$.  A \emph{$2$-tangle} is a pair $\left(B,t\right)$ such that the endpoints of the arcs go to a specific set of $4$ points on the surface of $B$, so that $2$ of the points lie on the upper hemisphere and the other $2$ points lie on the lower hemisphere with respect to the height function, and so that the interiors of the arcs are embedded in the interior of this ball. The two points in the upper hemisphere are referred to as the Northwest (NW) and Northeast (NE) corners respectively, and the two points in the lower hemisphere are referred to as the Southwest (SW) and Southeast (SE) corners respectively.  
\end{defn}

Rational tangles are a special class of $2$-tangles.

\begin{defn}[Rational tangle] 
A $2$-tangle is \emph{rational} if there exists an orientation preserving homeomorphism of pairs: 
$g:(B,t)\rightarrow (D^{2} \times I, \{x,y\} \times I)$ 
\end{defn}

It is easiest to visualize a tangle or a knot by its \emph{diagram}.  A diagram of a tangle or a knot can be viewed as a $4$-valent planar graph $G$, with over-under crossing information at each vertex.  For our knots in section \ref{sec:Constructions}, we will need to consider the number of \textit{twist regions}.  A twist region of a knot or tangle is a maximal string of bigons arranged from end to end in its diagram.  A single crossing adjacent to no bigons is also a twist region.  A diagram of a tangle is also embedded inside a great circle on which the endpoints of the tangle lie.  For a tangle $\left(B,t\right)$, assume that $B$ is the standard unit ball embedded in $\mathbb{R}^{3}$.  When we form the tangle diagram, $B$ projects to the standard unit disc $D$ in $\mathbb{R}^{2}$, with the upper hemisphere of B projecting to $D \cap \left\{ \left( x,y \right) \in \mathbb{R}^{2} : y \geq 0 \right\}$ and the lower hemisphere of B projecting to $D \cap \left\{ \left(x,y\right) \in \mathbb{R}^{2} : y \leq 0 \right\}$. Figure $1$ illustrates diagrams of two different $2$-tangles.  The left-hand side of Figure $1$ exhibits the diagram of a non-rational tangle inside its projection disc with its endpoints labelled.  The right-hand side shows a rational tangle with three twist regions.  Intuitively, we can differentiate between rational and non-rational $2$-tangles by the fact that rational tangles can be untwisted to the $0$-tangle, while non-rational tangles cannot.  

\begin{figure}[h]
\includegraphics[scale=0.65]{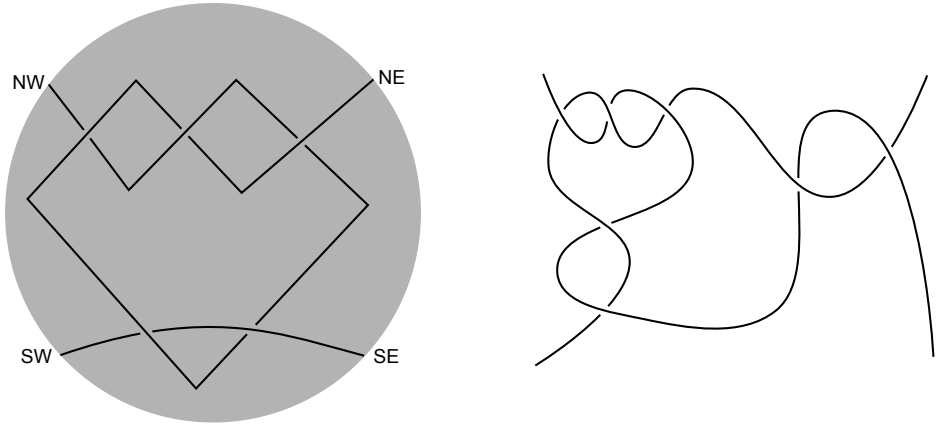}
\caption{The diagram of a non-rational $2$-tangle on the left and the diagram of a rational $2$-tangle on the right.}
\label{figure1}
\end{figure}

We say that two rational tangles, $T = \left(B,t\right) $ and $S = \left(B,s\right)$, are \emph{equivalent}, denoted $T$ $\sim$ $S$, if there exists an orientation-preserving homeomorphism $h:(B,t)\rightarrow (B,s)$ that is the identity map on the boundary.  Equivalently, $T$ $\sim$ $S$ if and only if they admit diagrams with identical configurations of their four endpoints on the boundary of the projection disc, and they differ by a finite number of Reidemeister moves, which take place in the interior of the disc.  A rational tangle is associated canonically with a unique reduced rational number or $\infty$, called \emph{the fraction of the tangle}.   The following theorem shows that rational tangles can be completely classified by their fractions.   

\begin{thm} (Conway, 1970)
\label{thm:Conway}
Two rational tangles are equivalent if and only if they have the same fraction.  
\end{thm}

For a proof, see \cite{KL}.  \\

We can give a topological meaning to this fraction.  For a given rational tangle $T = \left(B,t\right)$, there exists a compression disc embedded in $B$ and disjoint from the two arcs inside of $B$.  The fraction associated to $T$ corresponds to the slope of this compression disc.

For our purposes, we shall only be dealing with \emph{elementary rational tangles}, and even more specifically, \emph{vertical tangles}. We will combine these vertical tangles through the tangle operation of addition.

\begin{defn}[Elementary rational tangles] 
The \emph{integer tangles}, denoted by $\left[n\right]$, are made of $n$ horizontal twists, $n \in \mathbb{Z}$.
The \emph{vertical tangles}, denoted by $\frac{1}{\left[n\right]}$, are made of $n$ vertical twist, $n \in \mathbb{Z}$.  
Together, the integer tangles and vertical tangles make up the \emph{elementary  rational tangles}.  
\end{defn}

\begin{figure}[h]
\includegraphics[scale=0.80]{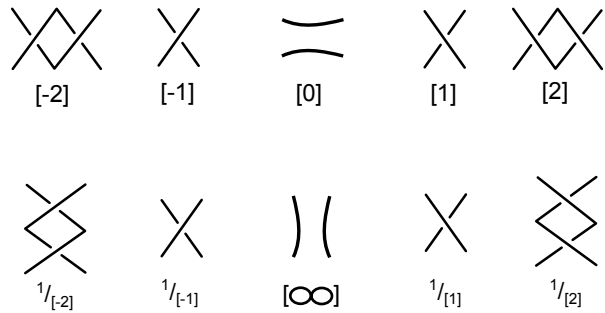}
\caption{Some of the elementary rational tangles.  Integer tangles are depicted in the top row and vertical tangles are depicted in the bottom row.}
\label{figure2}
\end{figure}

The tangle operation of addition is depicted below.  This operation is well-defined for $2$-tangles.  Note that, the sum of two rational tangles is a rational tangle if and only if one of the tangles is integral.  Thus, if we add two vertical tangles, the resulting tangle is non-rational.  This fact will be important for the proof of Theorem \ref{thm:knotvol}.

\begin{figure}[h]
\includegraphics[scale=0.40]{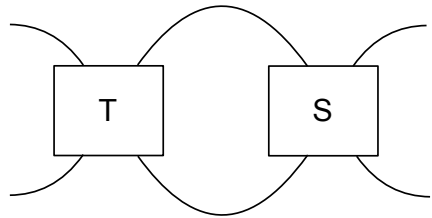}
\caption{The addition of two $2$-tangles.}
\label{figure3}
\end{figure}

Since by Theorem \ref{thm:Conway} we know that a rational tangle is completely determined by its fraction, we will often refer to a rational tangle by its associated fraction.

\subsection{Montesinos links} 
In this paper, we shall construct links by adding together a finite number of rational tangles.  Such links are called \emph{Montesinos links}.

\begin{figure}[h]
\includegraphics[scale=0.75]{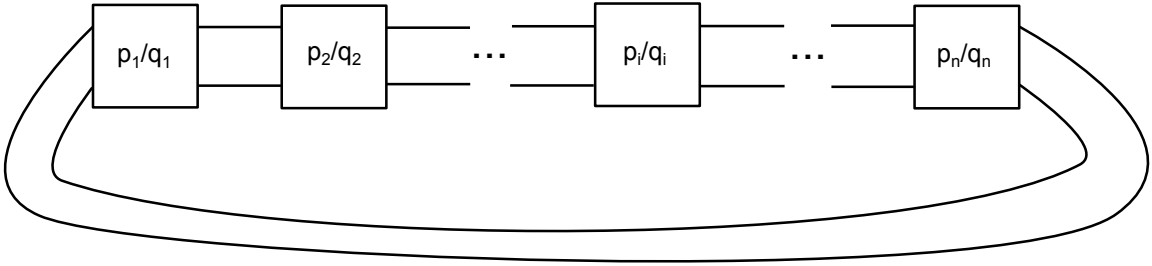}
\caption{A Montesinos link.  Each box denotes a rational tangle, and these tangles are connected to one another by tangle addition.}
\label{figure4}
\end{figure}

\begin{defn}[Montesinos link]
A \emph{Montesinos link}, denoted $K\left( \frac{p_{1}}{q_{1}}, \frac{p_{2}}{q_{2}},\ldots, \frac{p_{n}}{q_{n}} \right)$, is defined to be the link constructed by connecting $n$ rational tangles in a cyclic fashion, reading clockwise, with the $i^{th}$-tangle associated with the fraction $\frac{p_{i}}{q_{i}}$.  See figure \ref{figure4}.
\end{defn}

 We shall assume that $q_{i} \neq 1$ for each $i$.  If some $q_{i}=1$, then the $i^{th}$-tangle can be incorporated in an adjacent tangle, giving another Montesinos link:

\begin{center}
$K\left( \frac{p_{1}}{q_{1}},\ldots, \frac{p_{i-1}}{q_{i-1}}, \frac{p_{i}}{q_{i}},\frac{p_{i+1}}{q_{i+1}},\ldots, \frac{p_{n}}{q_{n}} \right)$ $=$ $K\left( \frac{p_{1}}{q_{1}},\ldots, \frac{p_{i-1}}{q_{i-1}}, \frac{p_{i}q_{i+1}+p_{i+1}}{q_{i+1}},\ldots, \frac{p_{n}}{q_{n}} \right)$
\end{center} 

The following classification of Montesinos links was originally proved by Bonahon in 1979 \cite{Bo}, and another proof was given by Boileau and Siebenmann in 1980 \cite{BS}.  A proof similar to the one done by Boileau and Siebenmann can be found in \cite[Theorem $12.29$]{BZ}. Also, Bedient gave a classification of a special class of Montesinos knots in \cite{Be}.

\begin{thm} \cite{Bo}
\label{thm:Bo}
The Montesinos links $K\left( \frac{p_{1}}{q_{1}}, \frac{p_{2}}{q_{2}},\ldots, \frac{p_{n}}{q_{n}} \right)$ and $K'\left( \frac{p'_{1}}{q'_{1}}, \frac{p'_{2}}{q'_{2}},\ldots, \frac{p'_{n}}{q'_{n}} \right)$ with $n \geq 3$ and $\sum_{j=1}^{n}\frac{1}{q_{j}} \leq n-2$, are classified by the ordered set of fractions $\left(\frac{p_{1}}{q_{1}} \, \mod \, 1, \ldots, \frac{p_{n}}{q_{n}} \, \mod \, 1\right)$ up to the action of the dihedral group generated by cyclic permutations and reversal of order, together with the rational number $\sum_{j=1}^{n} \frac{p_{j}}{q_{j}}$.
\end{thm}

\section{Mutations along Conway spheres}
\label{sec:Csphere}
In this section, we shall recall some results about mutations along Conway spheres, which will be needed to construct hyperbolic knots whose complements have the same volume.  

\begin{defn}[Conway sphere] 
\label{def:CSphere}
Let $K\subset S^{3}$ be a knot.  A \emph{Conway sphere} for $K$ is an embedded $2$-sphere meeting $K$ transversally in four points, i.e., a $4$-punctured sphere in the knot complement.   
\end{defn}

\begin{figure}[ht]
\includegraphics[scale=0.50]{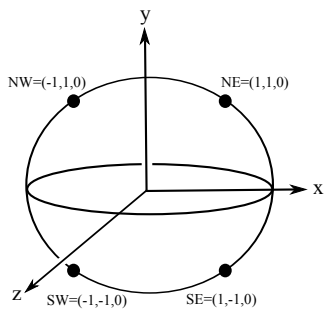}
\caption{A standard Conway sphere.}
\label{figure5}
\end{figure}

For a $4$-punctured sphere, denoted $S_{0,4}$, there exists a subgroup of its mapping class group isomorphic to $\mathbb{Z}/2\mathbb{Z} \times \mathbb{Z}/2\mathbb{Z} \subset Mod(S_{0,4})$, which is generated by the hyperelliptic (orientation preserving) involutions of $S_{0,4}$.  Note that the three non-trivial involutions are given by $180^{\circ}$ rotations about the $x$-axis, $y$-axis, and $z$-axis, respectively in figure \ref{figure5}.  For a given Conway sphere, $S$, we shall denote this subgroup by $\mathrm{Inv}(S)$.      

\begin{defn}[Mutation]
\label{def:Mutation}
A \emph{mutation} of a Conway sphere, $S$, in a $3$-manifold $M$ is the process of cutting $M$ along $S$ and then regluing $S$ back into $M$ by one of the nontrivial homeomorphisms of $\mathrm{Inv}(S)$.  If $K$ is a knot with a Conway sphere $S$, then cutting $(S^{3}, K)$ along $(S, S \cap K)$ and regluing by a mutation, $\mu$, yields a knot $K^{\mu} \subset S^{3}$.  We say that \emph{$K^{\mu}$ is obtained from $K$ by mutation along $S$}.
\end{defn}

Our examples in Section \ref{sec:KnotVol} are hyperbolic knots whose complements have the same volume.  The following theorem proved by Ruberman \cite{Ru} will be essential for these constructions.

\begin{thm} \cite{Ru}
\label{thm:Ru1}
Let $\mu$ be any mutation of a hyperbolic knot $K$ along an incompressible and $\partial$-incompressible Conway sphere.  Then $K^{\mu}$ is also hyperbolic, and the complements of $K$ and $K^{\mu}$ have the same volume.  
\end{thm}

The proof of this theorem relies on a careful analysis of the cusp geometry of hyperbolic $3$-manifolds and the embeddedness and uniqueness of least area surfaces in hyperbolic $3$-manifolds. Least area surfaces were first studied by Freedman--Hass--Scott; see \cite{FHS} and \cite{HS}.  These authors analyzed the compact case for least area surfaces, and Ruberman expanded this analysis to the non-compact case.  Also, a different proof of this theorem is given by Kuessner \cite{Ku}.

\section{Dehn surgery}
\label{sec:DSurgery}

Let $M$ be a $3$-manifold with torus boundary $\partial M$ and $s$ a \emph{slope} on $\partial M$.  The manifold obtained by gluing a solid torus $S^{1} \times D^{2}$ to $\partial M$ in such a way that the slope $s$ bounds a disc in the resulting manifold is called a \emph{Dehn surgery along s}.  In our case, manifolds with torus boundary will be knot complements.  

Given a hyperbolic $3$-manifold $M$ with a cusp corresponding to a torus boundary on $\partial M$, we choose a basis $\left\langle m,l \right\rangle$ for the fundamental group of the torus.  After this choice of basis, we can form the manifold $M\left(p,q\right)$ obtained by doing a $\left(p,q \right)$-Dehn surgery on the cusp, where $\left(p,q \right)$ is a coprime pair of integers.  By a \emph{$\left(p,q \right)$-Dehn surgery} we mean a Dehn surgery where we cut off the cusp and glue in a solid torus, mapping the boundary of the meridian disc to $s = pm+ql$.  

Thurston has shown that $M\left(p,q\right)$ is in fact a closed hyperbolic $3$-manifold for all $\left(p,q \right)$ near $\infty$; see \cite{Th}.  Neumann and Zagier expanded upon this analysis of Dehn surgery by giving an explicit formula relating the $\left(p,q \right)$ surgery coefficients to the geometry of the solid torus glued in during Dehn surgery:  

\begin{thm} \cite{NZ}
\label{thm:NZ}
Let $L$ be the length of the geodesic $\gamma$ on $M\left(p,q\right)$ which is in the homotopy class of the core of the solid torus $T$ added during the Dehn surgery. Let $Q\left( p,q \right) = \frac{\left( \mbox{length of pm+ql} \right)^{2}}{A}$, where $A$ denotes the area of the boundary torus.  Then 
\begin{center}
$L=2\pi Q\left( p,q \right)^{-1} + O\left(\frac{1}{p^{4}+q^{4}}\right)$
\end{center}
Thus, as the pair $\left(p,q \right) \rightarrow \infty$, $L \rightarrow 0$. 
\end{thm}

Our examples in Section \ref{sec:ClosedVol} are closed hyperbolic $3$-manifolds with the same volume, constructed by Dehn surgery on knot complements. For this, we make use of another result of Ruberman, which considers the effects of Dehn surgeries on the volume of a knot and its mutant.


\begin{defn}[Unlinked]
\label{def:Unlinked}
Let $K$ be a knot in $S^{3}$ admitting a Conway sphere $S$.  Observe that a specific choice of a mutation $\mu$ gives a pair of $S^{0}$'s on the knot such that each $S^{0}$ is preserved by $\mu$. We say that $\mu$ and $S$ are \textit{unlinked} if these $S^{0}$'s are unlinked on $K$.  
\end{defn}

\begin{thm} \cite{Ru}
\label{thm:Ru2}
Let $K$ be a hyperbolic knot in $S^{3}$ admitting an incompressible and $\partial$-incompressible Conway sphere $S$, and let $\mu$ be a mutation of S such that $\mu$ and $S$ are unlinked.  Then a $\left(p,q \right)$-Dehn surgery on $K$ and a $\left(p,q \right)$-Dehn surgery on $K^{\mu}$ give manifolds with the same volume.  
\end{thm}

\section{Constructing hyperbolic Montesinos knots}
\label{sec:Constructions}

\begin{figure}[ht]
\includegraphics[scale=0.75]{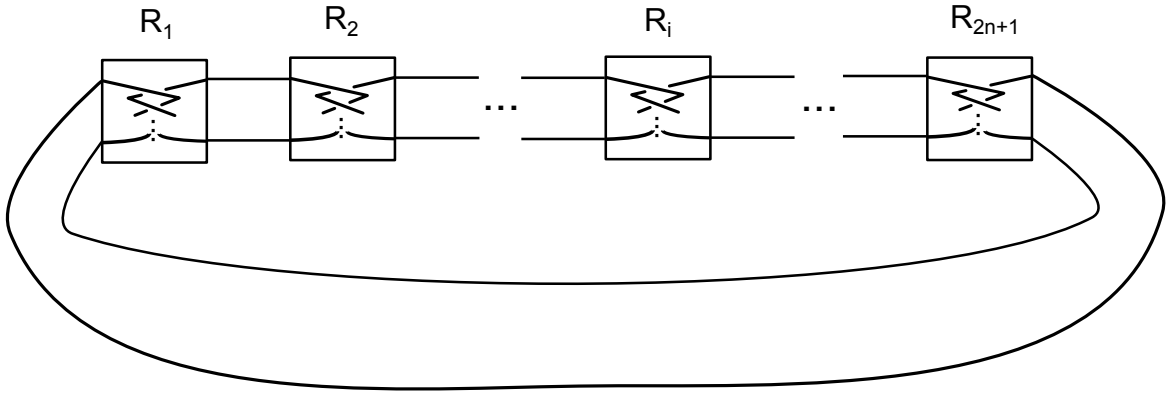}
\caption{A Montesinos link with $2n+1$ tangle regions.  Each tangle region $R_{i}$ contains a vertical tangle with $2i+5$ positive crossings.}
\label{figure6}
\end{figure}

Consider the Montesinos link, $K_{2n+1} = K\left( \frac{1}{7}, \frac{1}{9},\ldots, \frac{1}{2i+5}, \ldots, \frac{1}{4n+7} \right)$ depicted above.  Each $R_{i}$ in $K_{2n+1}$ represents a region in which the rational tangle $\frac{1}{ \left[ 2i+5 \right]}$ takes place, so this link is indeed a Montesinos link.  Note, for all $n \geq 1$, $K_{2n+1}$ is a link with at least two twist regions, each of which has at least $6$ crossings.  By the work of Futer and Purcell \cite[Theorem $1.4$]{FP}, this implies that each $K_{2n+1}$ is hyperbolic.  Next, it is easy to see that each rational tangle, $R_{i}$, is alternating, and since all our fractions are positive, i.e., all of our half twists performed in each $R_{i}$ region are positive, we have that each $K_{2n+1}$ is an alternating knot. Furthermore, we claim that each $K_{2n+1}$ is actually a knot and not a link with more than one component.  This follows from the fact that there are an odd number of vertical twists in each $R_{i}$ and there are an odd number of tangle regions.  An odd number of vertical twists in each $R_{i}$ results in one arc connecting the SE endpoint to the NW endpoint and the other arc connecting the NE endpoint to the SW endpoint in each tangle region.  The odd number of tangle regions ensures that $K_{2n+1}$ does not form two components when the $\left(2n+1\right)^{st}$ tangle connects back with the $1^{st}$ tangle.  We are focused on keeping our constructions knots because by the Gordon--Luecke Theorem \cite{GL}, we know that knots are determined by their complements.  

It is important to note that the above construction can be generalized.  In order to construct a hyperbolic, alternating Montesinos knot with $n$ twist regions ($n \geq 2$), $\widetilde{K}_{n}$, we just need our Montesinos link to meet a few parameters. Specifically, we need enough twisting in each twist region, an odd number of positive twists in each twist regions, and an odd number of twist regions.  In terms of a Montesinos knot's ordered set of fractions, this is equivalent to
\begin{center}
$\widetilde{K}_{n} = K \left( \frac{1}{t_{1}}, \frac{1}{t_{2}}, \ldots, \frac{1}{t_{2n+1}} \right)$, where $t_{i} > 6$ and $t_{i}$ is odd, for $1 \leq i \leq 2n+1$.   
\end{center}

\section{Hyperbolic knot complements that share the same volume}
\label{sec:KnotVol}
For the hyperbolic Montesinos knots $K_{2n+1}$ constructed in the previous section, consider the set of Conway spheres $ \left\{ S_{a} \right\}_{a=1}^{2n}$ where $S_{a}$ encloses only $R_{a}$ and $R_{a+1}$ on one side.  We can perform a mutation, $\sigma_{a}$, along $S_{a}$.  In particular, we have three possible choices of hyperelleptic involutions to use for a mutation along a Conway sphere.  For each $S_{a}$, we shall choose the mutation which interchanges NW with NE, and SW with SE, i.e., the rotation about the $y$-axis.  On one of our Montesinos knots, $K_{2n+1}$, such a mutation $\sigma_{a}$ interchanges the rational tangles $R_{a}$ and $R_{a+1}$.  In terms of this Montesinos knot's ordered set of fractions, such a mutation interchanges the $a^{th}$ fraction with the $\left(a+1\right)^{th}$.

Recall that $R_{a}$ and $R_{a+1}$ are vertical tangles, and so, they are rational, non-integral tangles.  We wish to apply Theorem \ref{thm:Ru1} to each of the mutations $\sigma_{a}$, $1\leq a \leq 2n$, which requires each Conway spheres $S_{a}$ to be both incompressible and $\partial$-incompressible.  This fact is established by the work of Oertel.

\begin{prop} \cite{Oe}
A Conway sphere $S$ is incompressible and $\partial$-incompressible in $S^{3} \setminus K$ if and only if it bounds at least two rational, non-integral tangles on each side.
\end{prop}

Putting these pieces together, we have the following proof of Theorem \ref{thm:knotvol}.

\begin{named}{Theorem \ref{thm:knotvol}}
There exists a sequence of volumes $\left\{v_{n}\right\}_{n=2}^{\infty}$ such that $N_{K}(v_{n}) \geq \frac{(2n)!}{2}$.

Furthermore, $\left(\frac{2n-1}{2}\right)v_{\mathrm{oct}} \leq v_{n}  \leq \left(4n+2\right)v_{\mathrm{oct}}$, where $v_{\mathrm{oct}} \left(\approx 3.6638\right)$ is the volume of a regular ideal octahedron.
\end{named}

\begin{proof}
For each $K_{2n+1}$, $n \geq 2$, consider the corresponding knot complement, $M_{2n+1} = S^{3} \setminus K_{2n+1}$.  By construction, we know that $M_{2n+1}$ is a hyperbolic $3$-manifold.  By Theorem \ref{thm:Ru1}, each $K_{2n+1}^{\sigma_{a}}$ is also a hyperbolic knot with $\mathrm{vol}(M_{2n+1}) = \mathrm{vol}(M_{2n+1}^{\sigma_{a}})$, for each mutation $\sigma_{a}$ performed along the Conway sphere $S_{a}$, $1 \leq a \leq 2n$.  Note that, any $S_{a}$ is incompressible and $\partial$-incompressible, since $n \geq 2$ implies that there are at least two rational, non-integral tangles on each side of $S_{a}$.  By repeatedly applying mutations, Theorem \ref{thm:Ru1} gives us that any finite word in our collection $ \left\{ \sigma_{a} \right\}_{a=1}^{2n}$ also preserves hyperbolicity and volume.  Since the set of adjacent transpositions generates the symmetric group, we have that $ \left\{ \sigma_{a} \right\}_{a=1}^{2n}$ generates $(2n+1)!$ possible compositions of mutations on $K_{2n+1}$ that preserve hyperbolicity and volume.  Next, by Theorem \ref{thm:Bo}, two of our Montesinos knots, $K_{2n+1}^{\sigma}$ and $K_{2n+1}^{\sigma'}$, where $\sigma$ and $\sigma'$ are finite words from our collection of $ \left\{ \sigma_{a} \right\}_{a=1}^{2n}$, are isomorphic if and only if they differ by a dihedral permutation of their rational tangles.  Since there are $2 \left( 2n+1 \right)= 4n+2$ dihedral permutations of our 2n+1 fractions, we have $\frac{(2n+1)!}{(4n+2)} = \frac{(2n)!}{2}$ different knots. By the Gordon--Luecke Theorem \cite{GL} we know that different knots have different knot complements, and so, we have the first result.  
 
As mentioned in section \ref{sec:Constructions}, each $K_{2n+1}$ is an alternating knot.  We can also easily check that each corresponding knot diagram is \emph{prime} and \emph{twist-reduced} by construction (see \cite{FP} for definitions of prime and twist-reduced).  This allows us to apply the work of Agol--Storm--Thurston \cite[Theorem $2.3$]{AST} to get lower bounds on volume for alternating hyperbolic knots, and we apply the work of Futer--Kalfagianni--Purcell \cite[Theorem $9.12$]{FKP} to get upper bounds on volume for Montesinos links. 
\end{proof}

\section{Closed hyperbolic $3$-manifolds that share the same volume}
\label{sec:ClosedVol}
In this section, we shall prove Theorem \ref{thm:closedvol} by constructing closed hyperbolic $3$-manifolds that share the same volume via Dehn surgery on our knot complements from the previous section.
  
\begin{named}{Theorem \ref{thm:closedvol}}
For each $n \geq 2$ there exists a convergent sequence $x_{n_{i}} \rightarrow v_{n}$ such that $N_{Cl}(x_{n_{i}}) \geq \frac{(2n)!}{2}$.  Furthermore, $x_{n_{i}} \leq \left(4n+2\right)v_{\mathrm{oct}}$, for all $i$. 
\end{named}
\begin{proof}  
By Theorem \ref{thm:knotvol}, we know that for $n \geq 2$, there exist $\frac{(2n)!}{2}$ different hyperbolic knots whose complements have the same volume, each of these knots coming from a finite number of mutations along Conway spheres in $K_{2n+1}$.  Let $\sigma$ and $\sigma'$ be any pair of distinct permutations generated by $ \left\{ \sigma_{a} \right\}_{a=1}^{2n}$.  Then $M^{\sigma}_{2n+1}$ and $M^{\sigma'}_{2n+1}$ are both hyperbolic knot complements with $\mathrm{vol}(M^{\sigma}_{2n+1})=\mathrm{vol}(M^{\sigma'}_{2n+1})$.  We wish to perform Dehn surgery on these manifolds and apply Theorem \ref{thm:Ru2}.

As noted in Theorem \ref{thm:Ru2}, in order to perform Dehn surgeries on a knot and its mutant and preserve volume, we need the corresponding Conway sphere and its mutation to be unlinked.  Since each $\sigma_{a}$ is a rotation about the $y$-axis of $S_{a}$, and each $S_{a}$ contains a tangle region in which one arc is connecting the SE endpoint to the NW endpoint and the other arc is connecting the NE endpoint to the SW endpoint, this unlinked condition will be met for all $S_{a}$ and $\sigma_{a}$, $1 \leq a \leq 2n$.  By Thurston's hyperbolic Dehn surgery theorem \cite{Th}, we know that for $(p,q)$ sufficiently large, $M^{\sigma}_{2n+1}\left(p,q\right)$ and $M^{\sigma'}_{2n+1}\left(p,q\right)$ are hyperbolic $3$-manifolds, where $(p,q)$ are the Dehn surgery coefficients.  Also, the cores of the solid tori filled in during these Dehn surgeries are isotopic to unique geodesics in their respective manifolds.  Thus, by Theorem \ref{thm:Ru2}, applying corresponding Dehn surgeries on $M^{\sigma}_{2n+1}$ and $M^{\sigma'}_{2n+1}$ yields $\mathrm{vol}(M^{\sigma}_{2n+1}\left(p,q\right)) = \mathrm{vol}(M^{\sigma'}_{2n+1}\left(p,q\right))$.  In order to complete the first part of this proof, we need to show that we can choose Dehn surgery coefficients $(p,q)$ so that $M^{\sigma}_{2n+1}\left(p,q\right)$ and $M^{\sigma'}_{2n+1}\left(p,q\right)$ are not homeomorphic hyperbolic $3$-manifolds, whenever $M^{\sigma}_{2n+1}$ and $M^{\sigma'}_{2n+1}$ are not homeomorphic knot complements.

Suppose not: that is, $M^{\sigma}_{2n+1}\left(p,q\right)$ is homeomorphic to $M^{\sigma'}_{2n+1}\left(p,q\right)$.  By Mostow rigidity \cite{BP}, we know there exists an isometry $f: M^{\sigma}_{2n+1}\left(p,q\right) \rightarrow  M^{\sigma'}_{2n+1}\left(p,q\right)$.  Since isometries map shortest geodesics to shortest geodesics, we will choose our $(p,q)$ so that the core geodesics, $\gamma$ and $\gamma'$, of the solid tori filled in during Dehn surgery correspond to the \textit{unique} shortest geodesics in their respective manifolds.  We are able to find such $(p,q)$ since by Theorem \ref{thm:NZ}, as $(p,q) \rightarrow (\infty,\infty)$, the lengths of the core geodesics approach $0$, while the other lengths stabilize.  Note that we are using the fact that the geodesic length spectrum of a hyperbolic $3$-manifold is discrete.  Thus, for our $(p,q)$ sufficiently large, the isometry $f: M^{\sigma}_{2n+1}\left(p,q\right) \rightarrow  M^{\sigma'}_{2n+1}\left(p,q\right)$ restricts to a homeomorphism $f': M^{\sigma}_{2n+1}\left(p,q\right) \setminus N(\gamma) \rightarrow M^{\sigma'}_{2n+1}\left(p,q\right) \setminus N(\gamma')$.  Fix $\epsilon > 0$ so that $\gamma$ and $\gamma'$ are the only geodesics in the $\epsilon$-thin parts of their respective manifolds. Then the previous statement would imply that $M^{\sigma}_{2n+1}\left(p,q\right)$ and $M^{\sigma'}_{2n+1}\left(p,q\right)$ have homeomorphic thick parts.  But these thick parts are homeomorphic to $M^{\sigma}_{2n+1}$ and $M^{\sigma'}_{2n+1}$, respectively. This contradicts the fact that $M^{\sigma}_{2n+1}$ and $M^{\sigma'}_{2n+1}$ are not homeomorphic.

Since for each $n \geq 2$, we have an infinite number of choices for $\left(p_{i},q_{i}\right)$-surgery coefficients large enough to meet our requirements, we get a sequence of volumes $\left\{x_{n_{i}}\right\}$ such that $N_{Cl}(x_{n_{i}}) \geq \frac{(2n)!}{2}$.  We obtain $x_{n_{i}} \rightarrow v_{n}$ from the fact that the filled manifolds $M^{\sigma}_{2n+1}\left(p_{i},q_{i}\right)$ converge to $M^{\sigma}_{2n+1}$.

The upper bound on volume follows from the upper bound for our hyperbolic knots in Theorem \ref{thm:knotvol}, and the fact that Dehn surgery strictly decreases volume \cite{Th}.
\end{proof}

\section{Conclusions and some open questions}
\label{sec:Final}

We are now ready to prove Theorem \ref{thm:mainthm}.  

\begin{named}{Theorem \ref{thm:mainthm}}
There exist sequences of volumes, $v_{n} , x_{n} \rightarrow \infty$ such that 
\begin{center}
$N_{K}(v_{n}) \geq \left(v_{n} \right)^{\left(\frac{v_{n}}{8} \right)}$ and $N_{Cl}(x_{n}) \geq \left(x_{n} \right)^{\left(\frac{x_{n}}{8} \right)}$
\end{center}
for all $n \gg 0$.
\end{named}

\begin{proof}
From Theorems \ref{thm:knotvol} and \ref{thm:closedvol}, we know that there exist sequences of volumes $\left\{v_{n}\right\}_{n=2}^{\infty}$ and $\left\{x_{n}\right\}_{n=2}^{\infty}$ with $v_{n} = \mathrm{vol}(M_{2n+1})$, $x_{n} = \mathrm{vol}(M_{2n+1}\left(p,q\right))$, such that \\ $v_{n}, x_{n} \leq \left( 4n+2 \right) v_{\mathrm{oct}}$ and $N_{K}(v_{n}), N_{Cl}(x_{n}) \geq \frac{(2n)!}{2}$ for all $n \geq 2$.  Using the estimate $n! \geq \sqrt{2\pi n} \left(\frac{n}{e} \right)^{n}$, which holds for all $n \geq 1$, and doing a little algebra yields:
\begin{center}
$N_{K}(v_{n}) \geq \frac{\left(2n\right)!}{2} > \left(\frac{2n}{e} \right)^{2n} \geq \left( \frac{v_{n}}{2ev_{\mathrm{oct}}}-\frac{1}{e} \right)^{ \frac{v_{n}}{2v_{\mathrm{oct}}} -1} > \left(\frac{v_{n}}{20} - \frac{1}{e} \right)^{\frac{v_{n}}{7.5} - 1} \geq \frac{v_{n}}{c}^{\frac{v_{n}}{7.5}}$, 
\end{center}
for some constant $c > 0$.
Thus, for $n \gg 0$, $N_{K}(v_{n}) \geq \left(v_{n}\right)^{\frac{v_{n}}{8}}$.  The same inequality holds for $N_{Cl}(x_{n})$.  
\end{proof} 

Theorem \ref{thm:mainthm} raises a number of questions.  To begin with, this theorem only deals with the classes of knot complements and closed manifolds with respect to $N(v)$.  So, can we find sequences $v_{n} \rightarrow \infty$ such that the growth rate of other special classes of hyperbolic $3$-manifolds with volume $v_{n}$ is exceptionally fast?  As mentioned in the introduction, for a few other special classes lower bounds and upper bounds on growth rates of $N_{v}$ have been established.  Frigerio--Martelli--Petronio \cite{FMP} have shown that the number of hyperbolic $3$-manifolds with totally geodesic boundary with volume $v$, grows at least as fast as $v^{cv}$ for some $c > 0$, Belolipetsky--Gelander--Lubotzky--Shalev \cite{BGLS} have shown that for $v \gg 0$, the number of arithmetic hyperbolic $3$-manifolds with volume $v$ is approximately $v^{c'v}$ for some $c' > 0$, and Lubotzky and Thurston \cite{BGLM}, \cite{Ca} have shown that the number of non-arithmetic hyperbolic $3$-manifolds grows at least factorially fast with $v$.  All of these examples, including our own, establish lower bounds for growth rates of $N(v)$ of the form $N_{SP}(v_{n}) \geq v_{n}^{cv_{n}}$, where SP denotes a special class of hyperbolic $3$-manifolds.  However, there has yet to be established a faster than factorial growth rate for $N(v)$ with respect to $v$ for any classes of hyperbolic $3$-manifolds.  This raises the question: does there exist a sequence of volumes $v_{n} \rightarrow \infty$ such that $N(v_{n})$ grows faster than factorially with $v$?    

We actually suspect that the answer is no.  From J\o rgensen's Theorem (found in \cite{Th}), we know that for any fixed volume $v > 0$, there are only finitely many homeomorphism types of the $\epsilon$-thick part of $M$, $M_{>\epsilon}$, where $\mathrm{vol}(M) \leq v$.  Consider the proof of J\o rgensen's Theorem and let us see how this can give us estimates on the number of $M_{>\epsilon}$ that can exist for $\mathrm{vol}(M_{>\epsilon}) \leq v$.  This technique is described by Burger--Gelander--Lubotzky--Mozes \cite{BGLM} who carefully consider possible finite covers of $M_{>\epsilon}$ by open balls.  In order to bound the number of thick parts $\left\{ M_{>\epsilon} \right\}$ bounded by a given volume $v$, it suffices to bound the number of possible fundamental groups for these thick parts.  This is done by associating to each such fundamental group a $2$-dimensional complex with the same fundamental group and counting these complexes.  The $2$-dimensional complex in question is the $2$-skeleton of the nerve of the open cover of $M_{>\epsilon}$.  As the authors point out, this work requires one to be careful in choosing the radii of the open balls in your cover sufficiently small with respect to the Margulis constant.  Then basic counting arguments show that there are at most $v^{cv}$ $2$-dimensional simpicial complexes for $\left\{ M_{>\epsilon} \right\}$  with volume at most $v$, for some constant $c$ which depends on the maximal degree of the vertices in the $1$-skeleton of our simplicial complexes, and so, there are at most $v^{cv}$ thick parts having volume $\leq v$.

In their paper, Burger--Gelander--Lubotzky--Mozes use the above argument to give similar bounds on the number of hyperbolic $n$-manifolds of volume $\leq v$ where $n \geq 4$.  They are easily able to do this since $\pi_{1} \left( M \right) = \pi_{1} \left( M_{> \epsilon} \right)$ for $n \geq 4$ by a simple application of Van Kampen's Theorem.  However, this quick argument does not work for hyperbolic $3$-manifolds.  To extend this upper bound on the growth rate of the number of general hyperbolic $3$-manifolds of volume $v$ (for some given $v$), we need to consider the effects of Dehn filling.  In particular, we need to know if the number of manifolds with a given volume is uniformly bounded amongst all hyperbolic Dehn fillings on any given cusped hyperbolic $3$-manifold.  This is equivalent to a difficult problem posed by Gromov in \cite{Gr}: is $N(v)$ locally bounded?

\bibliographystyle{hamsplain}
\bibliography{volumepaperbiblio}

\end{document}